\documentclass{amsart}
\usepackage{graphicx}
\usepackage{amscd}
\usepackage{amsmath}
\usepackage{amsfonts}
\usepackage{amssymb}
\usepackage{bbm}
\usepackage{setspace}
\usepackage{enumerate}         % better lists
\usepackage{fixme}
\usepackage{color}
\usepackage{url}
\usepackage{amsthm}
\usepackage{bm}
\usepackage{xy}
\usepackage{enumitem}

\theoremstyle{plain}
\newtheorem{theorem}{Theorem}[section]

\newtheorem{lemma}[theorem]{Lemma}

\newtheorem{proposition}[theorem]{Proposition}

\newtheorem{definition}[theorem]{Definition}

\theoremstyle{remark}
\newtheorem{remark}[theorem]{Remark}

\numberwithin{equation}{section}

\newcommand{\ind}{1\!\kern-1pt \mathrm{I}}
\newcommand{\rsto}{]\!\kern-1.8pt ]}
\newcommand{\lsto}{[\!\kern-1.7pt [}

\numberwithin{equation}{section}

\renewcommand{\rho}{\varrho}

\begin{document}
\title[A proof variant of FTAP]{A convergence result for the Emery topology and a variant of the proof of the fundamental theorem of asset pricing}
\begin{abstract}
We show that \emph{No unbounded profit with bounded risk} (NUPBR) implies  \emph{predictable uniform tightness} (P-UT), a boundedness property in the Emery topology which has been introduced by C.~Stricker~\cite{S:85}. Combining this insight with well known results from J.~M\'emin and L.~S\l{}ominski~\cite{MS:91} leads to a short variant of the proof of the fundamental theorem of asset pricing initially proved by F.~Delbaen and W.~Schachermayer~\cite{DS:94}. The results are formulated in the general setting of admissible portfolio wealth processes as laid down by Y.~Kabanov in~\cite{kab:97}.
\end{abstract}

\thanks{The second author gratefully acknowledges the support from
ETH-foundation. We thank Freddy Delbaen, Aleksandr Gushchin, Yuri Kabanov, Costas Kardaras, Irene Klein, Chong Liu and Walter Schachermayer for fruitful discussions on the topic.}

\author{Christa Cuchiero and Josef Teichmann}
\address{Vienna university of Technology, Wiedner Hauptstrasse 8--10, A-1040 Vienna and ETH Z\"urich, R\"amistrasse 101, CH-8092 Z\"urich}
\maketitle

\section{Introduction}

The single most important result of mathematical finance is the \emph{Fundamental Theorem of Asset Pricing} (FTAP): it establishes under a fairly weak assumption on a set $ \mathcal{X} $ of admissible portfolio wealth processes, a property called \emph{No Free Lunch with Vanishing Risk} (NFLVR), the existence of an equivalent separating measure $ Q \sim P $. This rather technical sounding assertion is the correct and sharp mathematical formulation of the vague ``meta-theorem'' stating that no arbitrage is \emph{essentially} equivalent to the existence of an equivalent martingale measure and has thus tremendous consequences: first, models can be easily characterized to satisfy (NFLVR) by simply checking whether such a separating measure $ Q \sim P $ exists. Second, the statement of FTAP is, mathematically speaking, the characterization of typical elements of a polar cone, which in turn allows to look at optimization problems from a dual point of view. Third, by simple economic arguments, separating measures $ Q \sim P $ lead to pricing structures for general payoffs. We believe that is worth searching for a simplification of proofs of FTAP: it might shed new light on these fundamentals of mathematical finance and and it makes the fundamentals more accessible to, e.g., students or practitioners. Additionally we believe that a simpler proof makes extensions, for instance towards robust finance or large financial markets, more reachable, see Remark \ref{robust_finance}.

For a detailed account of the long history of FTAP we refer to the overview article~\cite{S:10} by W. Schachermayer and the monograph~\cite{delsch:06} by F.~Delbaen and W.~Schachermayer. Let us here only briefly state (without any claim to completeness) the main milestones. The history of FTAP traces back to the work of F.~Black and M.~Scholes~\cite{BS:73} and R.~Merton in 1973. Indeed, their formula was the starting point for a deep investigation between the relation of pricing by no arbitrage considerations and pricing by taking ``risk neutral'' expectations (with respect to a martingale measure). In the late 1970s and early 1980s major advances in 
establishing a precise mathematical connection between those notions and proving first versions of FTAP in different settings were achieved by S.~Ross~\cite{R:78}, M.~Harrison, D.~Kreps and S.~Pliska~\cite{HK:79, HP:81, K:81}. These seminal papers have been generalized and further developed in many directions, in particular a first complete proof of FTAP in finite discrete time was given by R.~Dalang, A.~Morton and W.~Willinger~\cite{DMW:90} extending the Harrison-Pliska result~\cite{HP:81}. In continuous time, C.~Stricker~\cite{S:90} combined the result of D.~Kreps with a theorem by J.A. Yan~\cite{Y:80}, which is now known under the name Kreps-Yan theorem and which states the equivalence between \emph{No free lunch} (NFL) and the existence of an equivalent separating measure (see Theorem~\ref{krepsyan}).
The remaining major challenge was to replace the difficult interpretable and strong condition of (NFL) (involving closures in the weak-$*$-topology in $L^{\infty}$, see Section~\ref{sec:setting_kabanov} for the precise definition) by an economically convincing concept which only slightly strengthens the intuitive notion of absence of arbitrage. It turns out that the concept of (NFLVR) introduced in~\cite{DS:94} is precisely the right minimal and economically meaningful requirement which still allows to conclude the existence of an equivalent separating measure. Also the concept of \emph{No free lunch with bounded risk}, as applied in~\cite{Sch:94} in the discrete infinite time horizon case, would serve this purpose, however, (NFLVR) is weaker. The first complete proof of this ground-breaking result was presented by F.~Delbaen and W.~Schachermayer in~\cite{DS:94}, when the \emph{set of admissible portfolio wealth processes  $ \mathcal{X} $} is given by stochastic integrals $ (H \bullet S) $, for admissible integrands $H$ with respect to a finite dimensional locally bounded semimartingale $S$. Their beautiful and impressive proof builds on deep insights and is in some parts quite tricky, of which no essential simplification has been obtained since then. Finally in the accurate and sharply focused paper~\cite{kab:97}, Y.~Kabanov introduced, inspired also by~\cite{del:92}, an abstract setting of admissible portfolio wealth processes (see Definition~\ref{setting}): Y.~Kabanov's insight was that the proof of~\cite{DS:94} transfers almost literally to this more abstract setup and that there is a more condensed way of describing the proof. In particular, for obtaining the existence of an equivalent separating measure\footnote{In order to show that the separating measure is a local martingale measure, the local boundedness cannot be dropped.},  no local boundedness assumption on $S$ is necessary, which is also claimed in~\cite[Theorem 4.1]{DS:98}. The results of this article are formulated in the setting of~\cite{kab:97}, since we believe that it allows best for a simple presentation of the main findings.

The proof of FTAP is built upon an impressive series of technical lemmas (see Section \ref{sec:proof_ftap_kabanov}). By a slight change of strategy of proof we formulate two theorems (summarized in Theorem~\ref{essential_theorem} below) which allow for the same conclusion as the technical lemmas and which have a meaning on their own. In particular an interesting result about convergence in the Emery topology is proved. Let us be more precise on this issue: let $(X^n)$ be a sequence of admissible wealth processes where the final outcomes $X^n_1$ converge in probability. Under which circumstances can we conclude that actually $ X^n \to X $ in the (quite strong) Emery topology? If the processes were martingales and convergence in probability of final outcomes is replaced by convergence in $L^1$, then we could conclude by a very useful inequality due to Burkholder (see, e.g.,~\cite[Theorem 47, p.50]{mey:72}, or Section \ref{nupbr-to-put}) that for all predictable $ H $ with $ {\|H\|}_\infty \leq 1 $, and all $ a > 0 $
\[
aP\big[ {|(H \bullet (X^n-X))|_1}^* \geq a \big] \leq 18 {\|X^n_1-X_1\|}_1,  
\]
where $|X|_1^*=\sup_{0 \leq t\leq 1} |X_t|$. In~\cite{DS:94} and~\cite{kab:97} a series of lemmas is proved at this point to achieve -- with an additional maximality assumption -- the desired convergence in Emery topology under (NFLVR). We asked ourselves, whether it is possible to replace those lemmas by mathematical results with a meaning outside the proof of FTAP and with more financial interpretation. The main idea is to split the actual proof of convergence in the Emery topology into two steps: first, we show that a sequence of admissible portfolio wealth processes satisfying \emph{No Unbounded Profit with Bounded Risk} (NUPBR) also satisfies a boundedness property in the Emery topology. This boundedness property is known in the literature and is called \emph{predictable uniform tightness property} (P-UT) (see, e.g., ~\cite[VI.6.1]{JS:03}). Second, combining well-known results under (P-UT) together with the above mentioned maximality assumption (introduced in~\cite{DS:94}) we show that indeed convergence in the Emery topology follows. 

In~\cite{DS:94}, respectively~\cite{kab:97}, it is proved that (NFLVR) is equivalent to (NUPBR) together with the classical \emph{No Arbitrage} (NA) condition. This innocent looking insight has important consequences, since it allows to investigate the effects of (NUPBR) and (NA) separately: for instance in our proof the first main result is that already (NUPBR) implies (P-UT). We provide two proofs for this result: a more technical one mimicking some crucial arguments of~\cite{DS:94}, or~\cite{kab:97}, respectively, or a more financially inspired way using the existence of supermartingale deflators.  We do prefer the second one, since the existence of supermartingale deflators is equivalent to (NUPBR) (under fairly general assumptions, see~\cite{kar-GSD:13}) and very convincing from the point of view of applications: indeed in the realm of stochastic portfolio theory, see the excellent overview article~\cite{karfer:09}, it turned out to be very fruitful to analyze the consequences of (NUPBR) itself (in a setting of non-negative portfolio wealth processes), in particular to prove the crucial existence of supermartingale deflators, see, e.g.~\cite{KK:07},~\cite{rok:10} and~\cite{kar-GSD:13}. Having a supermartingale deflator at hand the proof of the (P-UT) property is indeed achieved via Burkholder's inequality, which restores the simplicity of the above martingale argument in the general situation.

The article is structured as follows: in Section \ref{sec:setting_kabanov} the setting of Y.~Kabanov's version of FTAP is introduced. In Section \ref{sec:proof_ftap_kabanov} the main steps of the proof in~\cite{kab:97} are presented, with precise references to the original proof in~\cite{DS:94}. In Section \ref{nupbr-to-put} the first main result of our work is derived, namely that (NUPBR) (or, in other terminology (NA1)) implies the (P-UT) property. In Section \ref{emery-convergence} we prove that (P-UT) together with a certain maximality assumption leads to convergence in the Emery topology, which -- by an easy argument -- concludes the proof of FTAP.

\section{Setting}\label{sec:setting_kabanov}

We work in the abstract setting of Y.~Kabanov as laid down in~\cite{kab:97}. We shall use for convenience of the reader~\cite{kab:97} as main reference, but we point out the corresponding results in~\cite{DS:94}. Let $\mathbb{S}$ be the space of semimartingales $X$ defined on a finite interval $[0,1]$ and starting from zero. The space $\mathbb{S}$ is equipped with the Emery topology defined by the metric
\[
d_E(X_1,X_2) := \sup_{K \in b \mathcal{E}, \, {\| K\|}_{\infty} \leq 1} E \big[ {|(K \bullet (X_1-X_2))|}^*_1 \wedge 1 \big] \, ,
\]
where $|X|_1^*=\sup_{t\leq 1}|X_t|$,  $b \mathcal{E}$ denotes the set of simple predictable strategies, that is, $K$ is of the form
\[
K=\sum_{i=0}^nK_i1_{]\tau_{i},\tau_{i+1}]} \, ,
\] 
with $n\in\mathbb{N}$, stopping times $0 = \tau_0 \leq \tau_1 \leq \dots \leq \tau_n \leq \tau_{n+1}=1$ and $K_i$ are $\mathcal{F}_{\tau_{i}}$-measurable random variables. The space of semimartingales is a complete metric space with the Emery topology, which follows essentially from the Bichteler-Dellacherie Theorem, see~\cite{eme:79}. Notice that M.~Emery defines the metric via the supremum over all bounded predictable processes (see~\cite{eme:79}) and not only over all simple predictable processes. However, as shown in, e.g.,~\cite{mem:80}, this leads to equivalent metrics.

Pathwise uniform convergence in probability is metrized by
\[
E[{|X-Y|}^*_1 \wedge 1 ] = d(X,Y) \, ,
\]
which makes the space of c\`adl\`ag processes a complete metric space. Obviously uniform convergence in probability is a weaker topology than the Emery topology.

The following definition is taken from~\cite{kab:97}:
\begin{definition}\label{setting}
Consider a convex set $\mathcal{X}_1 \subset \mathbb{S} $ of semimartingales
\begin{itemize}
\item starting at $0$, 
\item bounded from below by $-1$, 
\item being closed in the Emery topology, and 
\item satisfying the following \emph{concatenation property}: for all bounded, predictable strategies $ H,G \geq 0 $, $ X,Y \in \mathcal{X}_1 $ with $ HG=0 $ and $ Z = (H \bullet X) + (G \bullet Y) \geq -1 $, it holds that $ Z \in \mathcal{X}_1 $. 
\end{itemize}
We denote by $\mathcal{X}$ the set $ \mathcal{X} = \cup_{\lambda >0} \lambda \mathcal{X}_1 $  and call its elements \emph{admissible portfolio wealth processes}. The elements of $ \mathcal{X}_1 $ are called $1$-admissible wealth processes. We denote by $K_0$, respectively $K_0^1$ the evaluations of elements of $\mathcal{X}$, respectively $\mathcal{X}_1$, at terminal time $T=1$. 
\end{definition}

\begin{remark}\label{rem:settingDS}
Let $S$ be a $d$-dimensional semimartingale. Then the set of all stochastic integrals $ (\phi \bullet S) $, where $ \phi$ is $S$-integrable such that there exists a uniform bound from below $ (\phi \bullet S) \geq - \lambda $, for some $ \lambda \geq 0 $, is a set of admissible wealth processes generated by
\[
\mathcal{X}_1: =\big \{(\phi \bullet S) \, | \; \phi \text{ is } S\text{-integrable and} \, (\phi \bullet S) \geq -1 \, \big \}.
\]
The crucial properties of this set (axiomatized in the previous definition) are closedness in the Emery topology, which is a consequence of J.~M\'emin's theorem (see~\cite{mem:80}), convexity and the concatenation property, which are both just facts of stochastic integration theory. This is exactly the setting considered in~\cite{DS:98} (and also in~\cite{DS:94} with the additional requirement that $S$ is locally bounded).
\end{remark}

\begin{remark}
In line with the previous remark the set $\mathcal{X}$ of semimartingales introduced in Definition~\ref{setting} should be interpreted as \emph{discounted values} of portfolio wealth processes.
\end{remark}

\begin{remark}
We point out that the setting of Definition~\ref{setting} is structurally more general than the one of Remark~\ref{rem:settingDS}. An example provided in~\cite[Example 1.3]{Kar:12} shows that even in the one-period case under trading constraints a gap between num\'eraires and ``maximal'' portfolios appears. This cannot happen in the setting of Remark~\ref{rem:settingDS} as shown in~\cite{DS:95}.
\end{remark}

Let us introduce several notions of absence of arbitrage, for which we define the following convex cones:
\begin{align}\label{eq:cones}
C_0:=K_0 - L^0_{\geq 0}, \quad  C:=(K_0 - L^0_{\geq 0})\cap L^\infty.
\end{align}
\begin{description}
\item[(NA)] The set $\mathcal{X} $ is said to satisfy \emph{No Arbitrage} if 
\[
(K_0 - L^0_{\geq 0}) \cap L_{\geq 0}^0 = C_0 \cap L_{\geq 0}^0=\{0\},
\] 
which can be easily shown to be equivalent to 
\[
((K_0 - L^0_{\geq 0})\cap L^\infty) \cap L^{\infty}_{\geq 0} = C \cap L^{\infty}_{\geq 0} =\{0\}.
\]
\item[(NFLVR)] The set $\mathcal{X}$ is said to satisfy \emph{No free lunch with vanishing risk} if 
\[
\overline{C}\cap L_{\geq 0}^\infty=\{0\},
\]
where $\overline{C}$ denotes the norm closure in $L^{\infty}$.
\item[(NFL)] The set $\mathcal{X}$ is said to satisfy \emph{No free lunch} if 
\[
\overline{C}^* \cap L_{\geq 0}^\infty=\{0\},
\]
where $\overline{C}^*$ denotes the weak-$*$-closure in $L^{\infty}$.
\item[(NUPBR)] The set $\mathcal{X}_1$ is said to satisfy \emph{No unbounded profit with bounded risk} if $K^1_0$ is a bounded subset of $L^0$.
\end{description}

\begin{remark}
\begin{enumerate}
\item Note that in~\cite{kab:97} the (NUPBR) condition is introduced under the name (BK).
\item (NFLVR) can easily be proved to be equivalent to (NA) and (NUPBR), i.e., (NFLVR) $\Leftrightarrow$ (NA) + (NUPBR) (see e.g.,~\cite[Lemma 2.2]{kab:97}, or~\cite[Corollary 3.8]{DS:94}).
\item (NFLVR) or even (NUPBR) are economically convincing minimal requirement for models, but only (NFL) allows to conclude relatively directly the existence of an equivalent separating measure, defined below.
\end{enumerate}
\end{remark}

\begin{definition}
The set $\mathcal{X}$ satisfies the (ESM) (equivalent separating measure) property if there exists an equivalent measure $Q \sim P$ such that $\mathbb{E}_{Q}[X_1] \leq 0$ for all $ X \in \mathcal{X}$.
\end{definition}

Under (NFL), the (ESM) property is a consequence of the Kreps-Yan Theorem (see \ref{krepsyan} below), which in turn follows from
Hahn-Banach's Theorem. 

Apparently we have
\[
\text{(NFL)} \Longrightarrow \text{(NFLVR)} \Longrightarrow \text{(NA)} \, ,
\]
but it is an astonishing and deep insight that under (NFLVR) it holds that $ C = \overline{C}^* $, i.e.~the cone $C$ is already weak-$*$-closed and (NFL) holds.

In the formulation of~\cite{kab:97} the fundamental theorem of asset pricing reads as follows:

\begin{theorem}
Under (NFLVR) the cone $C$ is weak $*$-closed, hence (NFL) holds, which is equivalent to (ESM). In other words: (NVLVR) $\Leftrightarrow$ (ESM).
\end{theorem}

\begin{proof}
See Theorem 1.1 of~\cite{DS:94} in the case of a $d$-dimensional locally bounded semimartingale $S$ and Theorem 1.2 of~\cite{kab:97} in Y.~Kabanov's abstract setting.
\end{proof}

\section{The classical proof of FTAP}\label{sec:proof_ftap_kabanov}

In this section we sum up the main steps of the proof of FTAP as presented in~\cite{kab:97} for the setting of admissible portfolio wealth processes and indicate the corresponding results of~\cite{DS:94}. The proof splits in two parts. First a series of conclusions are presented, which can be easily motivated with financial (trading) arguments. Second five lemmas follow, whose content is more technical and which are hard to prove. It is the purpose of the present article to replace the second series of lemmas by a different line of arguments.

The first series of conclusions is the following:
\begin{enumerate}
\item The convex cone $C$ defined in~\eqref{eq:cones} is closed with respect to the weak-$*$-topology if and only if $C_0$ is Fatou-closed, i.e.~for any sequence $(f_n)$ in $C_0$ uniformly bounded from below and converging almost surely to $ f $ it holds that $ f \in C_0$, see beginning of Section 3 in~\cite{kab:97}, and~\cite[Theorem 2.1]{DS:94} essentially tracing back to A.~Grothendieck.
\item Take now $ -1 \leq f_n \in C_0 $ converging almost surely to $ f $. Then we can find $ f_n \leq g_n = Y^n_1 $ with $Y^n \in \mathcal{X}$. 
\item By (NA) it follows that each $ Y^n \in \mathcal{X}_1 $.
\item By (NUPBR) it follows that there are forward-convex combinations $ \widetilde{Y^n} \in \operatorname{conv}(Y^n,Y^{n+1},\ldots) $ such that $ \widetilde{Y^n_1} \to \widetilde{h_0} \geq f $ almost surely.
\item This implies that 
the set $\widehat{K}^1_0 \cap \{g \in L_0 \, |\, g \geq f\}$, where $\widehat{K}^1_0$ denotes the closure of $K_0^1$ in $L^0$, is non-empty. Since it is also bounded by (NUPBR) and closed, a maximal element $h_0$ exists (see beginning of Section 3 in~\cite{kab:97} or~\cite[Lemma 4.3]{DS:94}).
Since $h_0 \in \widehat{K}^1_0$, we can find a sequence of semimartingales $ X^n \in \mathcal{X}_1 $ such that $ X^n_1 \to h_0 $ almost surely and $ h_0 $ is maximal above $ f $ with this property. 
\item \label{convergence} The previously constructed ``maximal'' sequence of semimartingales $ X^n \in \mathcal{X}_1 $ converges pathwise uniformly in probability, i.e. $ {|X^n - X|}^*_1 \to 0 $ in probability, to some c\`adl\`ag process $X$ (see~\cite[Lemma 3.2]{kab:97} or~\cite[Lemma 4.5]{DS:94}).
\end{enumerate}

Since it is of crucial importance we devote a proper definition to maximality as in Section 3 of~\cite{kab:97} or before Lemma 4.3 in~\cite{DS:94}:
\begin{definition}\label{def:max}
An element $h_0 \in \widehat{K}^1_0$ (where $\widehat{K}^1_0$ denotes the closure of elements of $K^1_0$ which dominate $f$) is called maximal if it is maximal with respect to the pointwise (partial) ordering in $L^0$.
\end{definition}

It is now the goal to show that the sequence $(X^n)$ constructed in~\ref{convergence} above converges to $X$ in the Emery topology, an apparently much stronger statement. From this it follows that $h_0=\lim_{n \to \infty} X^n_1=X_1 \in K_0^1$, since \emph{$\mathcal{X}_1$ is closed in the Emery topology}. This it turn implies that $f \in C_0$, which finishes the proof by step (i) above.

Convergence in the Emery topology can be shown with respect to any equivalent measure $ Q \sim P $, since this notion of convergence only depends on the equivalence class of probability measures. By the basic convergence result \ref{convergence} we know that $ \xi := \sup_n {|X^n|}^*_1 \in L^0 $ (after passing to a subsequence). We can therefore find a measure $Q \sim P $ (take, e.g., $ dQ/dP = c \exp(-\xi) $) such that $ X^n \in L^2(Q) $, hence we can continue the analysis with $L^2$-methods, in order to prove Emery convergence with respect to $Q$.

Now the series of more technical lemmas starts: assume (NUPBR), take a sequence of (special) semimartingales $X^n=A^n+M^n$ whose sup-processes $|X^n|^*_1$ are uniformly bounded in $L^2(Q)$.
\begin{enumerate}
\item First key lemma: the sequence $ {|M^n|}_1^* $ is bounded in $L^0$ (see~\cite[Lemma 2.5]{kab:97} or~\cite[Lemma 4.7]{DS:94}). 
\item Second key lemma: define $ \tau^n_c := \inf\{ t \, | \; {|M^n|}_t^* > c \} $ for some $ c > 0 $, $ X^n_c := (1_{[\tau^n_c,\infty[} \bullet X^n) $, then for every $ \epsilon > 0 $ there is $ c_0 > 0 $ such that for all
\[
\widetilde{X} \in \cup_{c \geq c_0} \operatorname{conv}(X^1_c,\ldots,X^n_c,\ldots)
\]
it holds that $ Q[{|\widetilde{M}|}_1^* > \epsilon] \leq \epsilon $ (see~\cite[Lemma 2.6]{kab:97} or~\cite[Lemma 4.8]{DS:94}).
\item Third key lemma: for every $ \delta > 0 $ there is $ c_0 > 0 $ such that for all $ \widetilde{X} \in \cup_{c \geq c_0} \operatorname{conv}(X^1_c,\ldots,X^n_c,\ldots) $ it holds that $ d_E(\widetilde{M},0) \leq \delta $ (see~\cite[Lemma 2.7]{kab:97} or~\cite[Lemma 4.9]{DS:94}).
\item Fourth key lemma: there exist $ \widetilde{X}^n \in \operatorname{conv}(X_n,\ldots) $ such that $ \widetilde{M}^n $ converges in the Emery topology (see Lemma 2.8 in~\cite{kab:97} or Lemma 4.10 in~\cite{DS:94}).
\end{enumerate}

\begin{proposition}\label{emeryconvergence}
Let $\mathcal{X}_1$ satisfy (NUPBR). Let $ \widetilde{X}^n = \widetilde{M}^n+\widetilde{A}^n \in \mathcal{X}_1 $  be a sequence of special semimartingales, whose terminal values $X^n_1$ converge to a maximal element $ h_0 $ in probability such that $ \widetilde{M}^n$ converges in the Emery topology. Then $ \widetilde{A}^n$ converges in the Emery topology as well.
\end{proposition}

\begin{proof}
See~\cite[Lemma 3.3]{kab:97} or~\cite[Lemma 4.11]{DS:94}.
\end{proof}

As already argued above, this proposition together with key lemma (iv) implies that $ f \in C_0 $ yielding that $ C $ is in fact weak $*$-closed by step (i) above. Hence the assumptions of the Kreps-Yan Theorem \ref{krepsyan} are satisfied and we can conclude (ESM).

It is the purpose of this article to replace the arguments from key Lemma (i)--(iv) and Proposition \ref{emeryconvergence} by the following two theorems, which are sufficient to make the essential conclusions for the proof of FTAP. It is remarkable that no change of measure (which does a priori not have a financial interpretation) and no further passage to forward convex combinations is needed to achieve the result.

\begin{theorem}\label{essential_theorem}
Let $\mathcal{X}_1$ satisfy (NUPBR) and let $ (X^n)_{n \geq 0} \in \mathcal{X}_1 $ be a sequence of semimartingales. 
\begin{enumerate}
\item The sequence $(X^n)$ satisfies the (P-UT) property, i.e., 
\[
\lim_{a \to \infty} \sup_{H \in b\mathcal{E}, \|H\|_{\infty} \leq 1, \\ n \geq 0} P[|H \bullet X^n|_1 \geq a ] = 0 \, .
\] 
\item If the sequence $ (X^n) $ converges pathwise uniformly in probability to $X$ such that $X_1$ is a maximal element in $\widehat{K_0^1}$ (where $\widehat{K_0^1}$ denotes the closure of $K_0^1$ in $L^0$), then $ X^n \to X $ in the Emery topology.
\end{enumerate}
\end{theorem}

\begin{proof}
The proof is presented in Section \ref{nupbr-to-put} and \ref{emery-convergence}.
\end{proof}

\begin{remark}\label{robust_finance}
The presented proof has some implications on how to extend previous results towards new directions: for instance one could be tempted to replace Emery convergence by a weaker notion of convergence, where the supremum is taken only with respect to a smaller set of (more ``robust'') strategies, see for instance also~\cite{kar-GSD:13}. This also changes the (P-UT) property. With a properly adjusted notion of concatenation one will then be able to conclude a robust version of FTAP in general.
\end{remark}

\section{(NUPBR) implies (P-UT)}\label{nupbr-to-put}

Let us recall the (P-UT) property defined in~\cite[VI.6.1]{JS:03}. It was first considered by C.~Stricker~\cite{S:85} and taken up by A.~Jakubowski et al.~\cite{JMP:89} under the name (UT) (for ``uniforme tension''). The (P-UT) property has been studied by several authors, see in particular various papers by Kurtz and Protter~\cite{KP:91,KP:96}. 

\begin{definition}
We say that a sequence ${(X^n)}_{n \geq 0} $ of adapted c\`adl\`ag processes satisfies the (P-UT) property (predictably uniformly tight)
if for every $t >0$ the family of random variables $\{ (H \bullet X^n)_t\, :\,  H \in b\mathcal{E}, \|H\| \leq 1,  n \geq 0\}$ is bounded in $L^0$, that is,
\[
\lim_{a \to \infty} \sup_{H \in b\mathcal{E}, \|H\|_{\infty} \leq 1, \\ n \geq 0} P[|H \bullet X^n|_t \geq a ] = 0 \, .
\] 
\end{definition}

The (P-UT) property means that $(X^n)$ is bounded in the Emery topology.
From a mathematical finance point of view it can be seen as a new notion of absence of arbitrage by interpreting it 
as \emph{No unbounded profit or loss with simple predictable long or short positions in admissible portfolios.}

The heart of our considerations now consists in proving (NUPBR) implies (P-UT). From this conclusion it will be a relatively short way towards the existence of an equivalent separating measure: we provide two ways to show this result (Proposition~\ref{prop:NUPBR-PUT} and Appendix \ref{App:A}) . The first way builds upon the existence of a supermartingale deflator and the fact that sequences of supermartingales satisfy the (P-UT) property, which in turn easily translates to the original sequence of semimartingales. In Appendix~\ref{App:A} we provide a second approach by mimicking parts of the proof as presented in~\cite{kab:97} or~\cite{DS:94}.

\subsection{(NUPBR) implies (P-UT) -- an approach from mathematical finance} \label{sec:P-UT-MF}

It is well-known that the convex Emery-closed set $ 1 + \mathcal{X}_1 $ allows for a supermartingale deflator $ D $, see~\cite{rok:10},~\cite{kar-GSD:13},~\cite{kar-closure:13} and the recent preprint~\cite{imkper:14}, since a certain re-balancing property called \emph{fork-convexity} holds true.

\begin{definition}\label{def:forkconvex}
A set $\mathcal{Y} $ of semimartingales starting at $1$ is called \emph{fork-convex (with a strictly positive element)} if
\begin{itemize}
\item it contains a strictly positive semimartingale,
\item it is convex,
\item and for all times $ 0 \leq t \leq 1 $, for all $ A \in \mathcal{F}_t $ and for two wealth processes $ X,\tilde{X} \in \mathcal{Y} $, where $ \tilde{X} \in \mathcal{Y} $ is strictly positive, it holds that the portfolio,  obtained by re-investing the wealth of $X$ at time $ t $  into $\tilde{X} $ if $A$ happens, also belongs to $ \mathcal{Y} $. Formally speaking the portfolio defined through
\begin{itemize}
\item For $ 0 \leq s \leq t$ or $ \omega \notin A $: $ Y_s(\omega) :=X_s(\omega) $,
\item For $ t \leq s \leq T$ and $ \omega \in A$: $Y_s (\omega):= \frac{X_t(\omega)}{\tilde{X}_t(\omega)} \tilde{X}_s(\omega)$,
\end{itemize}
satisfies $ Y \in \mathcal{Y} $.
\end{itemize}
\end{definition}

\begin{remark}
Notice that the concatenation property of Definition~\ref{setting} is stronger than fork-convexity. The difference is comparable to the difference between simple predictable and predictable strategies.
\end{remark}

\begin{definition}
A positive c\`adl\`ag adapted process $ D $ is called supermartingale deflator for $ 1 + \mathcal{X}_1 $ if $D_0 \leq 1$ and $ D(1+X) $ is a supermartigale for all $ X \in \mathcal{X}_1$.
\end{definition}

\begin{theorem}\label{th:supermartingaledeflator}
Let $ \mathcal{X}_1  $ be a set of semimartingales satisfying Definition \ref{setting} and (NUPBR), then there exists a supermartingale deflator $ D $ for $ 1+\mathcal{X}_1$. 
\end{theorem}

\begin{remark}
Notice that a supermartingale deflator is essentially a process related to the first order condition of some utility optimization problem on $ 1 + K^1_0 $, see~\cite{rok:10},~\cite{kar-GSD:13} and~\cite{imkper:14}. One can also see this differently in the case where $ \mathcal{X}_1 $ is the set of integrals $ (H \bullet S) $ with respect to \emph{one} one-dimensional semimartingale $S$. In this case (and also under convex constraints on the trading strategies) the equivalence of (NUPBR) and the existence of a supermartingale deflator has been proved in~\cite[Proposition 4.12]{KK:07}. Looking at the characteristics of $S$ it is apparent that without a formal candidate for a density process via general Girsanov transforms, one gets ``immediate arbitrages'', i.e., not even (NUPBR) holds. This formal candidate, however, is necessarily a supermartingale deflator, since all involved processes are local martingales bounded from below.
\end{remark}

\begin{proof}
See, e.g.,~\cite[Lemma 2.3]{kar-closure:13}. To apply this result we have to prove that $ 1 + \mathcal{X}_1 $ is fork-convex. Indeed, fix $ 0 \leq t \leq 1 $ and $ A \in \mathcal{F}_t$, 
and choose two bounded predictable strategies, namely
\[
H := 1_{]0,t]}+1_{A^c}1_{]t,1]}
\]
and
\[
G^n := \left(\frac{X_t}{\tilde{X}_t}  \wedge n \right)\, 1_{A} 1_{]t,1]} \, ,
\]
which satisfy $ H,G^n \geq 0 $ and $ HG^n=0$ and
\begin{align*}
& (H \bullet X)_s(\omega) + (G^n \bullet \tilde{X} )_s(\omega) \\
= &\begin{cases}
X_s(\omega) -1 & \geq -1, \text{ if $ \omega \notin A$ or $0 \leq s \leq t$}\, , \\
X_t(\omega)-1 + \left(\frac{X_t(\omega)}{\tilde{X}_t(\omega)} \wedge n\right) \, (\tilde{X}_s(\omega) - \tilde{X}_t(\omega)) & \geq -1, \text{ else} \, . \end{cases}
\end{align*} 
Then it follows by the concatenation property of Definition \ref{setting} that 
\[
 (H \bullet X) + (G^n \bullet \tilde{X} ) \in \mathcal{X}_1 
\] 
for every $ n \geq 0 $. Since $ \mathcal{X}_1 $ is Emery closed by assumption, also the limit for $ n \to \infty $ lies in $ \mathcal{X}_1 $, but this limit equals $ Y - 1 $.
\end{proof}

The subsequent lemma shows that a sequence of non-negative  supermartingales with bounded initial values satisfies the (P-UT) property.

\begin{lemma}\label{lem:P-UTsuper}
Let $\mathcal{Z}$ be a set of non-negative supermartingales such that $Z_0 \leq K$ for all $ Z \in \mathcal{Z}$ and some $K > 0$. Then all sequences $(Z^n)_{n \geq 0}$ in $ \mathcal{Z}$ satisfy the (P-UT) property.
\end{lemma}
\begin{proof}
This is precisely the following elementary inequality by Burkholder, see, e.g.,~\cite[Theorem 47, p.50]{mey:72}. For every non-negative supermartingale $S$ and every process $H \in b\mathcal{E}$ with $\|H\|_{\infty} \leq 1$ it holds that
\[
a P[ {|(H \bullet S)|}^*_1 \geq a] \leq 9 E[|S_0|] 
\]
for all $ a \geq 0 $. Applying this inequality to $Z^n$ and letting $a \to \infty$ yields the (P-UT) property. Due to its importance for our arguments and as convenience for the reader, we provide a proof of Burkholder's inequality: fix $ a > 0 $, let $ S $ be a non-negative supermartingale and $ H $ bounded predictable with ${\|H\|}_{\infty} \leq 1 $, then $ Z:=S \wedge a $ is a supermartingale, too, and we have
\[
a P( {|(H \bullet S)|}^*_1 \geq a) \leq a P(|S|^*_1 \geq a) + a P({|(H \bullet Z)|}^*_1 \geq a) \, .
\]
Since $ Z $ is a supermartingale we obtain by the Doob-Meyer decomposition $ Z = M - A $ and $ (H \bullet Z) \leq (H \bullet M) + A $, which is a submartingale. Hence we can conclude by Doob's maximal inequalities for $p=2$ in case of the second term and $p=1$ (weak version) in case of the first term (see, e.g., Theorem 5.1 of~\cite{mey:72}) that
\[
a P( {|(H \bullet S)|}^*_1 \geq a) \leq E[S_0] +2 \frac{1}{a} E[{(H \bullet M)}^2_1 + A_1^2] \, .
\]
Ito's isometry allows to estimate the variance of the stochastic integral at time $1$ by $ E[M_1^2] $. Both quantities $M$ and $A$ of the Doob-Meyer decomposition may, however, be estimated through $ E[A_1^2] \leq E[M_1^2] \leq 2 a E[Z_0] $, since $Z$ is non-negative (so $A \leq M $ holds true) and $ Z \leq a $. This leads to the upper bound
\[
a P( {|(H \bullet S)|}^*_1 \geq a) \leq 9 E[S_0] \, .
\]
The estimate $ E[M_1^2] \leq 2 a E[Z_0] $ on the Doob-Meyer decomposition is well-known and follows from~\cite[Theorem 46]{mey:72}.
\end{proof}

\begin{remark}
The condition of bounded initial values can be dropped, if convergence in law of the supermartingales $(Z^n)$ on the Skorokhod space is assumed.
The corresponding statement can be found in~\cite[Proposition 3.2 c)]{JMP:89}.
\end{remark}

The following lemma is inspired by~\cite[Lemma 2.9]{kar-closure:13}, where convergence in the Emery topology of the integrator processes is replaced by the (P-UT) property (which corresponds to boundedness in the Emery topology) and uniform convergence in probability of the integrands by boundedness in $L^0$. The lemma is needed to establish the final proof that $(X^n)$ satisfies (P-UT).

\begin{lemma}\label{lem:P-UTint}
Let $(S^n)_{n \geq 0}$ be a sequence of semimartingales satisfying the (P-UT) property and let $(H^n)_{n \geq 0}$ be a sequence of adapted c\`adl\`ag processes
such that  $ (|H^n|^*_1)_{n \geq 0}$ is bounded in $L^0$. Then $(H^n_{-} \bullet S^n)_{n \geq 0}$ satisfies the (P-UT) property.
\end{lemma}
\begin{proof}
The assertion follows from~\cite[Corollary 6.20 b)]{JS:03} by noting that boundedness of $(|H^n|^*_1)_n$ in $L^0$ implies
boundedness of $( |H^n_{-}|^*_1)_n$ in $L^0$ since $|\Delta H^n|^*_1\leq 2|H^n|^*_1$. For the convenience of the readers we provide a proof here: 
consider the supremum process $ {|H^n|}^* $, which is bounded at time $ T =1 $ by assumption. Hence the stopping times
\[
\tau^n_k := \inf \{ s \leq 1 \, | \,  {|H^n|}^*_s \geq k \}
\]
are localizing as $ k \to \infty $ for any fixed $n \in \mathbb{N}$. Apparently $ {(H^n_{-} \bullet S^n)}^{\tau^n_k} $ satisfies the (P-UT) property, for every $k$, by assumption on $ (S^n) $, since $ | H^n_{t}| \leq k $ for $ t < \tau^n_k $. Now we can conclude that even  $(H^n_{-} \bullet S^n) $ satisfies (P-UT): indeed fix $ \varepsilon > 0 $, then we find a number $ k $ large enough such that $ P(\tau^n_k < \infty) \leq \varepsilon $ by assumption on $H^n$, for every $n$, and we find a constant $ C > 0 $ (not depending on $n$) by the (P-UT) property such that
\[
P[ |H \bullet S^n|_t \geq C ] \leq \varepsilon
\]
for $H\in b\mathcal{E}$ with $ \| H \|_{\infty} \leq 1 $, and every $n$. Hence
\[
P[ |(HH^n_{-} \bullet S^n)|_t \geq Ck ] \leq P[ |(HH^n_{-} \bullet S^n)|^{\tau^n_k}_t \geq Ck ] + P[\tau^n_k < \infty] \leq 2 \varepsilon  \, ,
\]
which concludes the proof.
\end{proof}

\begin{proposition}\label{prop:NUPBR-PUT}
Let $\mathcal{X}_1$ satisfy (NUPBR) and let $ (X^n)_{n \geq 0} \in \mathcal{X}_1 $ be \emph{any} sequence of semimartingales. Then $(X^n)$ satisfies the (P-UT) property.
\end{proposition}

\begin{remark}
The conclusion still holds true if the concatenation property of Definition~\ref{setting} is replaced by the weaker fork convexity of Definition~\ref{def:forkconvex}.
\end{remark}

\begin{proof}
Theorem~\ref{th:supermartingaledeflator} implies the existence of a supermartingale deflator and by Lemma~\ref{lem:P-UTsuper} the sequence of supermartingales $(Z^n):=(D(1+X^n))$ satisfies the (P-UT) property since $Z_0^n \leq 1$. By It\^o's integration by parts formula, we have
\[
1+X^n=\left(\frac{1}{D}\right)_{-} \bullet Z^n+Z^n_{-} \bullet \frac{1}{D}+ \left[\frac{1}{D}, Z^n\right].
\]
Applying Lemma~\ref{lem:P-UTint} twice, first with $H^n=\frac{1}{D}$ and $S^n=Z^n$ and second with $H^n=Z^n$ and $S^n=\frac{1}{D}$, we deduce that the first two terms satisfy the (P-UT) property. Note that the second application of Lemma~\ref{lem:P-UTint} is justified since the (P-UT) property  of $(Z^n)$ implies that $(|Z^n|^*_1)_n$ is bounded in $L^0$ (see e.g.~\cite[Lemma 1.3]{MS:91} or Proposition~\ref{prop:PUT-L0-bound} below).
The third term $[\frac{1}{D}, Z^n]$ also satisfies the (P-UT) property as
\[
\left[\frac{1}{D}, Z^n\right]_t=\frac{1}{2}\left(\left[\frac{1}{D}+ Z^n, \frac{1}{D}+ Z^n\right]_t-[ Z^n,Z^n]_t-\left[\frac{1}{D}, \frac{1}{D}\right]_t\right)
\]
and each of the terms on the right hand side is bounded in $L^0$ which is again a consequence of the (P-UT) property of $(Z^n)$ (see e.g.~\cite[Lemma 1.3]{MS:91} or Proposition~\ref{prop:PUT-L0-bound}).
\end{proof}

\section{A convergence result for semimartingales in the Emery topology}\label{emery-convergence}

This section is devoted to the proof of the second assertion of Theorem \ref{essential_theorem}. The strategy of the proof is to use the established (P-UT) property under (NUPBR) and to work with the decomposition of a semimartingale into big jumps and a special semimartingale. By a result of M\'emin and S\l{}ominski~\cite{MS:91}, the (P-UT) property allows to conclude Emery convergence for all except the finite variation part of the special semimartingale. For its convergence we need the maximality property.

At this point it is remarkable that already the (P-UT) property, which is easily established by the existence of supermartingale deflators, yields the essential part of the convergence result in the Emery topology.

\begin{theorem}\label{th:Emeryconv}
Let $\mathcal{X}_1$ satisfy (NUPBR) and let $ (X^n)_{n \geq 0} \in \mathcal{X}_1 $ be a sequence of semimartingales, which converges pathwise uniformly in probability to $X$ such that $X_1$ is a maximal element in $\widehat{K_0^1}$, where $\widehat{K_0^1}$ denotes the closure of $K_0^1$ in $L^0$. Then $ X^n \to X $ in the Emery topology.
\end{theorem}

\begin{proof}
Due to Proposition~\ref{prop:NUPBR-PUT}, (NUPBR) implies the (P-UT) property of $(X^n)$. Hence the theorem is a consequence of Proposition~\ref{prop:Memin} and Proposition~\ref{prop:kabanovmax} below.
\end{proof}

For a sequence of semimartingales $(X^n)_{n \geq 0}$ with $X^n_0=0$ and some $C >0$ let us consider the following decomposition
\begin{align}\label{eq:decomp}
X^n= B^{n,C}+M^{n,C}+\check{X}^{n,C},
\end{align}
where  $\check{X}^{n,C}= \sum_{s\leq t} \Delta X^n_s 1_{\{|\Delta X^n_s| > C\}}$, 
$ B^{n,C}$ is the predictable finite variation part and  $M^{n,C}$ the local martingale part of the canonical decomposition of the special semimartingale $X^n-\check{X}^{n,C}$.

The following proposition is a reformulation of~\cite[Proposition 1.10]{MS:91} and establishes the announced Emery convergence of the local martingale and the big jump part under (P-UT) for sequences of semimartingales  which converges pathwise uniformly in probability.

\begin{proposition}\label{prop:Memin}
Let $ (X^n)_{n \geq 0} $ be a sequence of semimartingales with $X^n_0=0$, which converges pathwise uniformly in probability to $X$.  Assume furthermore the (P-UT) property for this sequence and consider decompositions of form~\eqref{eq:decomp} for $(X^n)$ and $X$. Then there exists some $C >0$ such that $M^{n,C} \to M^C$ and $\check{X}^{n,C} \to \check{X}^{C}$ in the Emery topology and $B^{n,C}\to B^C$ pathwise uniformly in probability.
\end{proposition} 

\begin{proof}
Let us first note that $\Delta X^n \to \Delta X$ pathwise uniformly in probability. Since there exists some $C$ such that $P\left[\exists t \, |\, |\Delta X_t| =C\right]=0$ as stated in Lemma~\ref{lem:existC} below, it follows that
\[
 \sum_{s\leq 1} |\Delta X^n_s 1_{\{|\Delta X^n_s| > C\}}-\Delta X_s1_{\{|\Delta X^n_s| > C\}}| \to 0,
\]
implying convergence in variation and thus convergence in the Emery topology.
The remaining part follows from~\cite[Proposition 1.10]{MS:91}. This latter result is based on~\cite[Corollaire 1.9]{MS:91} (see also Theorem~\ref{th:intconv} below), from which convergence in probability of $[M^{n,C}-M^C, M^{n,C}-M^C]_1 \to 0$, and thus Emery convergence of $(M^{n,C})$, can be deduced.
\end{proof}

\begin{lemma}\label{lem:existC}
Let $X$ be a semimartingale. Then there exists some $C >0$ such that $P\left[\exists t \, |\, |\Delta X_t| =C\right]=0$.
\end{lemma}
\begin{proof}
The sets $\Omega_C := \left\{ \exists t \, |\, |\Delta X_t| =C \right\} $ are well-defined and measurable for $ C \geq 0 $. Furthermore their union is $ \Omega $, hence only countably many of these sets can have positive probability.
\end{proof}

\begin{remark}
In the case of a sequence of special semimartingales $(X^n)$ with canonical decomposition $X^n=A^n+N^n$, A.~Gushchin~\cite[Theorem 3]{Gus:96} proved that under
an additional condition, called (UI2), the conclusion of Proposition~\ref{prop:Memin} also holds for the parts of the canonical decomposition, namely that the finite variation part  $(A^n)$ converges uniformly in probability  and that the local martingale part $(N^n)$ converges in the Emery topology (as long as the (P-UT) property and uniform convergence in probability are satisfied). The condition (UI2) is for example implied if $(|\Delta X^n|^*_t)_{n\geq 0}$ is uniformly integrable for each $t \in [0,1]$.
By performing an equivalent measure change to $Q \sim P$ (as in the original proof of FTAP, see the paragraph below Definition~\ref{def:max})
such that $\sup_n |X^n|^*_1 \in L^2(Q)$, one could also deal with special $L^2$-semimartingales which satisfy in particular the above (UI2) condition, since $(|\Delta X^n|^*_t)_{n\geq 0}$ is uniformly integrable for each $t \in [0,1]$ due to $|\Delta X^n|^*_t\leq 2 | X^n|^*_t\leq 2\sup_n |X^n|^*_1$. An application of Gushchin's result would then yield Emery convergence of the local martingale part of the canonical special semimartingale decomposition, too.
\end{remark}

\begin{proposition}\label{prop:kabanovmax}
Let $\mathcal{X}_1$ satisfy (NUPBR) and let $ (X^n)_{n \geq 0} \in \mathcal{X}_1 $ be a sequence of semimartingales, which converges pathwise uniformly in probability to $X$ such that $X_1$ is a maximal element in $\widehat{K_0^1}$. (Here, $\widehat{K_0^1}$ denotes the closure of $K_0^1$ in $L^0$.) Consider the semimartingale decompositions of form~\eqref{eq:decomp} for $(X^n)$ and $X$ and
assume that $M^{n,C} \to M^C$ and $\check{X}^{n,C} \to \check{X}^{C}$ in the Emery topology. Then $B^{n,C}\to B^C$ in the Emery topology.
\end{proposition}

\begin{proof}
The proof goes along the lines of~\cite[Lemma 3.3]{kab:97} following~\cite[Lemma 4.11]{DS:94} by replacing $M^n$ by the slightly more involved expression $Y^n:=M^{n,C}+\check{X}^{n,C}$, which converges in the Emery topology by assumption. For the reader's convenience we however provide the full proof.
Let $r^n:=\frac{dB^{n,C}}{dB}$ where $B$ is a predictable increasing process dominating all $B^{n,C}$. Assume by contradiction that $(B^{n,C})$ does not converge in the Emery topology. Then it is not a Cauchy sequence and there exists $i_k, j_k \to \infty$ such that
\[
P\left[(|r^{i_k}-r^{j_k}| \bullet B)_1 > 2\gamma\right]\geq 2 \gamma >0.
\]
Let $\Gamma_k:=\{r^{i_k}\geq r^{j_k}\}$. Then we may conclude -- by assuming $i_k \wedge j_k \geq i_{k-1} \vee j_{k-1}$ and eventually interchanging $i_k$ and $j_k$ -- that
\[
P\left[(|r^{i_k}-r^{j_k}|1_{\Gamma_k} \bullet B)_1 > \gamma\right]\geq  \gamma >0.
\]
Take $\alpha_k \downarrow 0$ and define $\bar{X}^k:=1_{\Gamma_k}\bullet X^{i_k}+1_{(\Gamma_k)^c}\bullet X^{j_k}$ and 
\[
\sigma_k=\inf\{ t\, |\, (1_{\Gamma_k}\bullet Y^{i_k})_t+(1_{(\Gamma_k)^c}\bullet Y^{j_k})_t < Y_t^{i_k} \vee Y_t^{j_k}-\alpha_k\}.
\]
Note that $\bar{Y}^{k}-Y^{i_k}=1_{(\Gamma_k)^c}\bullet (Y^{j_k}-Y^{i_k})$ and $\bar{Y}^k-Y^{j_k}=1_{\Gamma_k}\bullet (Y^{i_k}-Y^{j_k})$
converges to $0$ in the Emery topology and thus also uniformly along path in probability. We  may therefore take $i_k$ and $j_k$ growing fast enough to ensure $P[\sigma_k <\infty] \to 0$. 
Set now $\widetilde{X}^k:=1_{[0,\sigma_k]}\bullet\bar{X}^k$.
By Lemma~\ref{lem:admissibility}, $\widetilde{X}^k\in (1+\alpha_k) \mathcal{X}_1$ and we have the following representation
\begin{align*}
\widetilde{X}^k_1&=(1_{\Gamma_k \cap [0, \sigma_k]}\bullet X^{i_k})_1+(1_{(\Gamma_k)^c \cap [0, \sigma_k]}\bullet X^{j_k})_1\\
&=X^{j_k}_{1\wedge \sigma_k}+(1_{\Gamma_k \cap [0, \sigma_k]}\bullet(X^{i_k}-X^{j_k}))_1\\
& =X^{j_k}_{1\wedge \sigma_k}+(1_{\Gamma_k \cap [0, \sigma_k]}\bullet (Y^{i_k}-Y^{j_k}))_1+\xi_k,
\end{align*}
where $\xi_k=(1_{\Gamma_k\cap [0, \sigma_k]}\bullet(B^{i_k,C}-B^{j_k,C}))_1=(1_{\Gamma_k}(r^{i_k}-r^{j_k})\bullet B)_{1\wedge \sigma_k}$.
Applying~\cite[Lemma A]{kab:97} to $\xi_k$, implies that forward convex combination of $\xi_k$ converge to a random variable $\eta\geq 0$ with $\eta \neq 0$. Denoting the maximal element to which $X_1^{j_k}$ converges by $h_0 \in \widehat{K}_0^1$, it follows - by the Emery convergence of $(Y^n)$ -- that forward convex combinations of $\widetilde{X}_1^k$ converge to $h_0 +\eta$. Since $\widetilde{X}^k\in (1+\alpha_k)\mathcal{X}_1$, this yields a contradiction to the maximality of $h_0 \in \widehat{K}_0^1$.
\end{proof}

\begin{lemma}\label{lem:admissibility}
Let $X^1, X^2 \in \mathcal{X}_1$ and $\alpha >0$. Consider the decomposition~\eqref{eq:decomp} for $X^1$ and $X^2$, i.e.,
\[
X^i=B^{i,C}+M^{i,C}+\check{X}^{i,C}=:B^{i,C}+Y^i, \quad  i=1,2.
\]
Let $B$ be a predictable increasing process dominating $B^{i,C}$. Set $r^i=\frac{dB^{i,C}}{dB}$ and 
\[
\sigma:=\inf\{t \, |\, (1_{\Gamma}\bullet Y^1)_t+(1_{\Gamma^c} \bullet Y^2)_t < Y^1_t \vee Y^2_t-\alpha\},
\]
where $\Gamma=\{r^1 \geq r^2\}$. Then 
\[
\widetilde{X}=1_{\Gamma \cap [0, \sigma]}\bullet X^1+1_{\Gamma^c \cap [0, \sigma]}\bullet X^2 \in (1+\alpha)\mathcal{X}_1.
\]
\end{lemma}

\begin{proof}
Notice that $1_{\Gamma} \bullet B^{1,C}+1_{\Gamma^c} \bullet B^{2,C} \geq B^{1,C} \vee B^{2,C}$. Thus on $[0,\sigma[$
\[
\widetilde{X}\geq B^{1,C} \vee B^{2,C}+Y^1 \vee Y^2-\alpha\geq (B^{1,C}+Y^1)\vee(B^{2,C}+Y^2)-\alpha=X^1 \vee X^2 -\alpha.
\]
At time $\sigma$ the jump of $\Delta\widetilde{X}$ is given by $\Delta\widetilde{X}_{\sigma}=1_{\Gamma }\Delta X_{\sigma}^1+1_{\Gamma^c }\Delta X_{\sigma}^2$ and hence $\widetilde{X}_{\sigma}\geq -1-\alpha$ because the left limit is at least $X^1 \vee X^2 -\alpha$.
\end{proof}

\appendix
\section{ The P-UT property}

This section is dedicated to state some results related to the (P-UT) property, which can be found in~\cite{MS:91},~\cite{JMP:89} and 
\cite[Chapter VI.6]{JS:03}. 

\begin{proposition}\label{prop:PUT-L0-bound}
Let $ (X^n)_{n \geq 0} $ be a sequence of semimartingales satisfying the (P-UT) property. Then 
\begin{enumerate}
\item 
$(|X^n|^*_1)_{n\geq 0}$ is bounded in $L^0$;
\item $([X^n, X^n]_1)_{n \geq 0}$ is bounded in $L^0$.
\end{enumerate}
\end{proposition}

\begin{proof}
See~\cite[Lemme 1.2]{JMP:89} or~\cite[Lemme 1.3]{MS:91}.
\end{proof}

The following proposition characterizes the (P-UT) property in terms of $L^0$-boundedness of certain parts in the semimartingale decomposition~\eqref{eq:decomp}.

\begin{proposition}
Let $ (X^n)_{n \geq 0} $ be a sequence of semimartingales and consider the semimartingale decomposition~\eqref{eq:decomp} for some $C >0$.
Then $ (X^n) $ satisfies the (P-UT) property if and only if
the following three conditions hold:
\begin{enumerate}
\item The sequence of total variations of $\check{X}^{n,C}$ denoted by $ {(\operatorname{TV}(\check{X}^{n,C})_1)}_{n \geq 0}$  is bounded in $L^0$.
\item The sequence $([M^{n,C}, M^{n,C}]_1)_{n \geq 0}$ is  bounded in $L^0$.
\item The sequence of total variations of $B^{n,C}$ denoted by $ {(\operatorname{TV}(B^{n,C})_1)}_{n \geq 0}$  is bounded in $L^0$.
\end{enumerate}
\end{proposition}

\begin{proof}
See~\cite[Th{\'e}or{\`e}me 1.4]{MS:91} or~\cite[Theorem VI.6.15]{JS:03}.
\end{proof}

The following theorem builds the basis of Proposition~\ref{prop:Memin}.

\begin{theorem}\label{th:intconv}
Let $(H^n)_{n \geq 0}$ be a sequence of adapted c\`adl\`ag processes and $(X^n)_{n \geq 0} $  a sequence of semimartingales satisfying the (P-UT) property. If $(H^n,X^n)$ converges pathwise uniformly in probability to $(H,X)$, then the stochastic integrals  $(H^n_{-} \bullet X^n)$ converge to  $(H_{-} \bullet X)$  pathwise uniformly in probability as well. In particular, $[X^n, X^n]_1 \to [X,X]_1$ in probability. 
\end{theorem}

\begin{proof}
See~\cite[Th{\'e}or{\`e}me 2.6]{JMP:89} or~\cite[Th{\'e}or{\`e}me 1.8, Corollaire 1.9]{MS:91}.
\end{proof}

\section{ (NUPBR) implies (P-UT) -- a direct approach} \label{App:A}

In this section we provide a direct proof of the (P-UT) property under the (NUPBR) condition. In contrast to Section~\ref{sec:P-UT-MF}, we here additionally need to assume pathwise uniform convergence in probability of the sequence $(X^n)$.

\begin{proposition}
Let $\mathcal{X}_1$ satisfy (NUPBR) and let $ (X^n)_{n \geq 0} \in \mathcal{X}_1$ be a sequence of semimartingales with $X^n_0=0$, which converges pathwise uniformly in probability to an adapted, c\`adl\`ag process $X$. 
\begin{enumerate}
\item Then for every $ C>0 $ there exists a decomposition $ X^n=M^n+B^n+\check{X}^n$ into a local martingale $M^n$, a predictable, finite variation process $B^n$ and a finite variation process $\check{X}^n$, for $ n \geq 0 $, such that jumps of $M^n$ and $B^n$ are bounded by $ 2C $ uniformly in $n$. 
\item  The sequence $ {(\operatorname{TV}(\check{X}^n)_1)}_{n \geq 0}$ of total variations of $\check{X}^n$ is bounded in $L^0$
and ${(\check{X}^n)}_{n \geq 0}$ satisfies the (P-UT) property.
\item The sequence $(|M^n|^*_1)_{n \geq 0}$ is bounded in $L^0$ and  $(M^n)_{n \geq 0}$ satisfies the (P-UT) property.
\item The sequence $ {(\operatorname{TV}(B^n)_1)}_{n \geq 0}$ of total variations of $B^n$ is  bounded in $L^0$ and $(B^n)_{n \geq 0}$ satisfies the (P-UT) property.
\item The sequence $(X^n)_{n \geq 0}$ satisfies the (P-UT) property.
\end{enumerate}
\end{proposition}

\begin{proof} 
Assertion (i) is simply a consequence of the semimartingale decomposition of form~\eqref{eq:decomp}, where we leave away the superscript $C$ for notational convenience.

Concerning (ii), we apply similar arguments as in~\cite[Corollaire 2.9]{JJ:81}. Define the c\`adl\`ag moduls for a c\`adl\`ag function $f$ by 
\[
w'(f, \delta)=\inf_{\Pi} \max_{i \leq r}\sup_{s,t \in [t_{i-1}, t_i)} |f(s)-f(t)|,
\]
where the infimum runs over all partitions $\Pi=\{0 =t_0 < t_1 \cdots < t_r=1\}$ with $\min_{i} (t_i-t_{i-1}) \geq \delta$.
Let $\varepsilon >0$. Since $(X^n)$ converges uniformly in probability to a c\`adl\`ag process $X$
 there exists some $\delta >0$ such that
\[
\sup_nP\left[w'(X^n,\delta) \geq C\right] \leq \sup_nP\left[w'(X^n-X,\delta)+w'(X,\delta) \geq C\right] \leq \frac{\varepsilon}{2}
\]
and some $b$ such that
\[
\sup_n P\left[|X^n|_1^* \geq b\right] \leq \frac{\varepsilon}{2}. 
\] 
If $w'(X,\delta) < C$, the definition of the c\`adl\`ag moduls implies that the number of jumps of $X$ in $[0,1]$, whose absolute value is greater than $C$, is at most $\lceil \frac{1}{\delta} \rceil$. Moreover the jumps of $X$ can be estimated by
\[
\sup_{t \leq 1} |\Delta X_t| \leq 2 |X|_1^*.
\]
Therefore 
\[
\sup_n P\left[\operatorname{TV}(\check{X}^n)_1 \geq \frac{2b}{\delta}\right] \leq \varepsilon, 
\] 
which implies the first claim. Moreover $(\check{X}^n)$ satisfies the (P-UT) property, since 
\[
\sup_{H \in b\mathcal{E}, \|H\|_{\infty} \leq 1, \\ n \geq 0} P[{|H \bullet \check{X}^n|}_t \geq a ]\leq 
\sup_{ n \geq 0} P[\operatorname{TV}(\check{X}^n)_1  \geq a ],
\] 
which converges due to the first assertion to $0$ as $a \to \infty$.
 
In order to prove the first part of assertion (iii), we follow the proof of~\cite[Lemma 2.5]{kab:97} (compare ~\cite[Lemma 4.7]{DS:94}), however without the assuming $L^2$-boundedness of $|X^n|^*_1$. For the reader's convenience we here repeat the main arguments and adjust it to our setting. Suppose that $(|M^n|^*_1)$ is unbounded in $L^0$. Then there exists a subsequence such that $P[|M^n|^*_1 > n^3] \geq 7\alpha$ for some $\alpha >0$. As $(X^n)$ converges uniformly in probability, there exists some $K \in \mathbb{R}_+$ such that
\[
P \left[|X^n|^*_1 > K\right ]\leq \alpha,
\]
in particular for $n$ large enough we have $P \left[|X^n|^*_1 > n\right ]\leq \alpha$. Let 
\[
\tau_n=\inf\{ t \geq 0\,| \, |M^n|^*_t \geq n^3 \textrm{ or } |X^n|^*_t \geq n\}
\] 
and $\widetilde{X}_n :=n^{-3} 1_{[0,\tau_n]} X^n$. We then have  $P[|\widetilde{M}^n|^*_1 \geq 1] \geq 6\alpha$ and $|\Delta \widetilde{M}^n|\leq 2C n^{-3}\leq n^{-1}$, since the jumps of $M^n$ are bounded by $2C$. Define stopping times $T^n_i$ inductively by 
\[
T^n_0=0, \quad T^n_{i+1}=\inf\{t \geq T^n_i \, |\, |\widetilde{M}^n_t -\widetilde{M}^n_{T_i}| \geq n^{-1} \}, \quad i \geq 0.
\]
By Lemma~\ref{lem:slicing} we thus have
\[
P[\widetilde{M}^n_{T^n_i}- \widetilde{M}^n_{T^n_{i-1}} < -\alpha n^{-1}] \geq \alpha  \quad \textrm{for all }i=1,\ldots, k_n:=\left\lfloor \frac{\alpha n}{2}\right\rfloor.
\]
Define $Z^n=X^n-\check{X}^n$. Since $|X^n|_1^*$ and $|\check{X}^n|_1^*$ are bounded in $L^0$, $|Z^n|_1^*$ is bounded in $L^0$ as well.
Thus we have for $n$ large enough and arbitrary stopping times $\sigma, \tau$,
\[
P \left[|\widetilde{Z}^n_{\sigma} -\widetilde{Z}^n_{\tau}|>\frac{\alpha n^{-1}}{2}\right]\leq P \left[2|Z^n|_1^* > \frac{\alpha n^{2}}{2}\right] \leq \frac{\alpha}{2}.
\]
Hence 
\[
P\left[\widetilde{B}^n_{T^n_i}- \widetilde{B}^n_{T^n_{i-1}} > \frac{\alpha n^{-1}}{2}\right] \geq \frac{\alpha}{2}. 
\]
Let $H^n=1_{[0,T^n_{k_n}] \cap \{r_n =1\}}$ where $r_n:=d \widetilde{B}^n/d|\widetilde{B}^n|$. Set $Y^n=H^n \bullet \widetilde{X}^n$. Then
$H^n \bullet \widetilde{B}^n$ and $H^n \bullet \widetilde{B}^n - \widetilde{B}^n$ are increasing on $[0, T^n_{k_n}]$ and since $1_{]T^n_{i-1},T^n_{i}]}H^n \bullet \widetilde{B}^n > \widetilde{B}^n_{T^n_i}- \widetilde{B}^n_{T^n_{i-1}}$, we have
\[
P\left[(1_{]T^n_{i-1},T^n_{i}]}H^n \bullet \widetilde{B}^n)_1 > \frac{\alpha n^{-1}}{2}\right] \geq \frac{\alpha}{2}. 
\]
By~\cite[Lemma B]{kab:97} it follows that
\[
P\left[(H^n \bullet \widetilde{B}^n)_1 > \frac{k_n \alpha^2 n^{-1}}{8}\right] \geq \frac{\alpha}{4},
\]
that is 
\begin{align}\label{eq:estimateB}
P\left[(H^n \bullet \widetilde{B}^n)_1 > 2 \beta \right] \geq 2 \beta
\end{align} 
for some $\beta >0$.
Since $X^n$ on $[0, \tau_n[$ is in $[-1, n]$, we have $\Delta X^n \geq -(n+1)$ on $[0, \tau_n]$. Thus 
\[
(\Delta Y^n)^- \leq  (\Delta \widetilde{X}^n)^- \leq \frac{(n+1)}{n^3}.
\]
Consider now $\|(H^n \bullet \widetilde{M}^n)_1\|_2^2$. As $|\widetilde{M}^n_{T^n_i}-\widetilde{M}^n_{T^n_{i-1}}| \leq \frac{2}{n} $ a.s., we can estimate
\begin{align*}
\|(H^n \bullet \widetilde{M}^n)_1\|_2^2 & = \left\|\left(\left(\sum_{i=1}^{k_n} 1_{]T^n_{i-1}, T^n_i]}\right)H^n \bullet \widetilde{M}^n\right)_1\right\|_2^2\\
&\leq \sum_{i=1}^{k_n} \|(1_{]T^n_{i-1}, T^n_i]}H^n \bullet \widetilde{M}^n)_1\|_2^2  
& \leq \frac{4 k_n}{n^2}.
\end{align*}
Since the local martingale $H^n \bullet \widetilde{M}^n$ is uniformly bounded (by $\alpha$), it is a true martingale and  Doob's inequality yields
\[
P[| H^n \bullet  \widetilde{M}^n|^*_1 \geq \gamma_n] \leq \gamma_n^{-2} \|(H^n \bullet \widetilde{M}^n)_1\|_2^2= \gamma_n^{-2}\frac{4 k_n}{n^2},
\]
which tends to $0$ for $\gamma_n^{-1}=o(\sqrt{n})$. Since $Y^n \geq  H^n \bullet  \widetilde{M}^n +  H^n \bullet  \widetilde{\check{X}^n}$, we have 
\begin{align*}
P\left[|(Y^n)^-|^*_1 \geq \gamma_n\right] &=P\left[\inf_{t \leq 1} Y_t^n \leq -\gamma_n\right]\\
&\leq
P\left[\inf_{t\leq 1}  ((H^n \bullet  \widetilde{M}^n)_t +  (H^n \bullet  \widetilde{\check{X}^n})_t) \leq -\gamma_n\right]\\
&\leq P\left[| H^n \bullet  \widetilde{M}^n +  H^n \bullet  \widetilde{\check{X}^n}|^*_1 \geq \gamma_n\right]\\
&\leq P\left[| H^n \bullet  \widetilde{M}^n |^*_1 + \frac{1}{n^3}|H^n \bullet  \check{X}^n|^*_1 \geq \gamma_n \right].
\end{align*}
By the above estimate and since $(\check{X}^n)$ satisfies the (P-UT) property as proved in step (ii), this tends to $0$. Due to~\eqref{eq:estimateB}, $P[Y^n_1 \geq \beta] \geq \beta$ and~\cite[Lemma 2.4]{kab:97} yields a contradiction to (NUPBR) and thus $(|M^n|_1^*)_{n \geq 0}$ is bounded in $L^0$.\\
Concerning the second assertion of (iii), it follows from~\cite[Proposition VI.6.13 (iv) $\Rightarrow$ (i)]{JS:03} that $(M^n)$ satisfies the (P-UT) property. Let us here repeat the arguments of this proof.
We have for all stopping times $T$
\[ 
 \mathbb{E}[[ M^n, M^n]_T]\leq \mathbb{E}[|(M^n)^2|^*_T].
\]
Since $\Delta(|(M^n)^2|^*_t)\leq C^2$,  Lenglart's domination property~\cite[Lemma I.3.32]{JS:03} implies
\begin{align*}
P\left[[ M^n, M^n]_1 \geq a  \right] &\leq \frac{1}{a}\left (\eta + \mathbb{E}[\sup_{t \leq 1}\Delta(|(M^n)^2|^*_t)]\right)+ P[|(M^n)^2|^*_1 \geq \eta]\\
&\leq\frac{1}{a}\left (\eta + C^2\right)+P[|(M^n)^2|^*_1 \geq \eta],
\end{align*}
which tends to $0$ for $a \to \infty$ and $\eta=\sqrt{a}$. Hence $[ M^n, M^n]_1$ is bounded in $L^0$. For $H^n \in b\mathcal{E}$ with $\|H^n\|_{\infty} \leq 1$, we have for all stopping times $T$
\[
\mathbb{E}[(H^n\bullet M^n)^2_T]=\mathbb{E}[(H^n)^2 \bullet [M^n, M^n]_T]
\]
and an  another application of Lenglart's domination property  yields
\begin{align*}
P\left[(H^n\bullet M^n)^2_1 \geq a  \right] &\leq \frac{1}{a}\left (\eta + \mathbb{E}[\sup_{t \leq 1}\Delta([M^n, M^n]_t)]\right)+ P[[ M^n, M^n]_1\geq \eta]\\
&\leq \frac{1}{a}\left (\eta + C^2\right)+P[[ M^n, M^n]_1\geq \eta],
\end{align*}
which tends again to $0$ for $a \to \infty$ and $\eta=\sqrt{a}$. Hence $(M^n)$ satisfies the (P-UT) property.

Let us now continue to prove statement (iv), namely that $ {(\operatorname{TV}(B^n)_1)}_{n \geq 0}$  is bounded in $L^0$. Suppose the claim is false. Decomposing $B^n$ into two increasing predictable processes such that $B^n = C^n -D^n$ and $\operatorname{TV}(B^n)=C^n+D^n$, this means that  $C^n$ or\slash and $D^n$ is unbounded in $L^0$. If $D^n$ is unbounded in $L^0$, we get that $C^n$ is unbounded in $L^0$ as well, since $X^n=M^n+C^n-D^n+\check{X}^n$ is bounded from below by $-1$ and $M^n$ and $\check{X}^n$ are bounded in $L^0$. Thus there exists a subsequence such that $P[C^n_1 > n^3] \geq 2\alpha$ for some $\alpha >0$ and as before we have for $n$ large enough  $P \left[|X^n|^*_1 > n\right ]\leq \alpha$. Let 
\[
\sigma_n=\inf\{ t \geq 0\,| \, C_t^n\geq n^3 \textrm{ or } |X^n|^*_t \geq n\}
\]
and $\widehat{X}_n :=n^{-3} 1_{[0,\sigma_n]} X^n$.  We then have  $P[\widehat{C}^n_1 \geq 1] \geq \alpha$. Note that
$\widehat{C}^n =(G^n \bullet \widehat{B}^n)_1$ where $G^n=1_{\{r_n =1\}}$ with $r_n:=d \widehat{B}^n/|\widehat{B}^n|$. Define $W^n=G^n \bullet \widehat{X}_n $. As $(M^n)$ and $(\check{X}^n)$ satisfy the (P-UT) property we obtain similarly as before
\begin{align*}
&P\left[W^n_1 >  \beta \right] \geq \beta, \quad \textrm{for some } \beta >0,\\
&(\Delta W^n)^- \leq \frac{(n+1)}{n^3},\\
&P\left[|(W^n)^-|^*_1 \geq \gamma_n\right] \to 0, \quad \textrm{for } \gamma_n \downarrow 0,
\end{align*}
which contradicts (NUPBR) due to~\cite[Lemma 2.4]{kab:97}  and thus  $ {(\operatorname{TV}(B^n)_1)}_{n \geq 0}$  is bounded in $L^0$.
Moreover $(B^n)$ satisfies the (P-UT) property, since
\[
\lim_{a \to \infty}\sup_{H \in b\mathcal{E}, \|H\|_{\infty} \leq 1, \\ n \geq 0} P[{|H \bullet B^n|}^*_1 \geq a ]\leq 
\lim_{a \to \infty}\sup_{ n \geq 0} P[\operatorname{TV}(B^n)_1  \geq a ] =0.
\] 
Assertion (v) follows from the fact that each of the sequences $(B^n)$, $(M^n)$ and $(\check{X}^n)$ in the decomposition of $(X^n)$ satisfies the (P-UT) property.
\end{proof}

The following lemma is an adaptation of~\cite[Appendix, Lemma D]{kab:97} where the requirement of an $\mathcal{H}^2$-martingale is replaced by the local martingale property together with uniformly bounded jumps.

\begin{lemma}\label{lem:slicing}
Let $\varepsilon > 0$ and let $N$ be a local martingale with $|\Delta N| \leq \varepsilon$. Define stopping times $T_i$ inductively by 
\[
T_0=0, \quad T_{i+1}=\inf\{t \geq T_i \, |\, |N_t -N_{T_i}| \geq \varepsilon \}, \quad i \geq 0.
\]
Assume that $P[|N|_1^* \geq 1 ] \geq 6 \alpha $ for some $0< \alpha <1$. Then 
\[
P[N_{T_i}- N_{T_{i-1}} < -\alpha \varepsilon] \geq \alpha \quad \textrm{for all }i=1,\ldots, k:=\left\lfloor \frac{\alpha}{2 \varepsilon}\right\rfloor.
\]
\end{lemma}

\begin{proof}
Let $f_i=N_{T_i}- N_{T_{i-1}}$ and $\Gamma=\{T_k \leq 1 \}$. 
As $|N_t -N_{T_{i-1}}| < \varepsilon $ a.s. for $T_{i-1}< t < T_i$, we have
\[
|1_{]T_{i-1}, T_i]} \bullet N|_1^*=\sup_{t \leq 1}\left|\int_0^t 1_{]T_{i-1}, T_i]} dN_s\right|\leq \sup_{T_{i-1} <t \leq 1} |N_{t\wedge T_i}- N_{T_{i-1}}| \leq |f_i| \leq 2 \varepsilon \quad \textrm{a.s.},
\]
where the last inequality follows from  the definition of the stopping times $T_i$ and the assumption on the jumps.
Using this, we have
\[
|N|^*_1 1_{\Gamma^c} = \left|\left(\sum_{i=1}^k 1_{]T_{i-1}, T_i]}\right) \bullet N\right|_1^*1_{\Gamma^c}\leq \sum_{i=1}^k |1_{]T_{i-1}, T_i]} \bullet N|_1^* \leq 2 k\varepsilon \leq \alpha <1, \quad \textrm{a.s.}
\]
Hence
\[
P[\Gamma] \geq P[|N|_1^* \geq 1 ] - \underbrace{P[\{|N|_1^* \geq 1\} \cap \Gamma^c ]}_{=0} \geq 6 \alpha.
\]
Let $L^i_t=N_{t \wedge T_i}-N_{t \wedge T_{i-1}}$. Since $(L^i_t)$ is a uniformly bounded (by $2 \varepsilon$) local martingale, it is a martingale and since
$L^i_{T_i}=f_i$, we have $\mathbb{E}[f_i]=0$. Therefore $\mathbb{E}[f_i^+]= \mathbb{E}[f_i^-]$ which implies  $\mathbb{E}[|f_i|]=2\mathbb{E}[f_i^-]$. As for all $i \leq k$,  $\{|f_i| \geq \varepsilon\} \supseteq \Gamma$, we have by Markov's inequality
\[
\mathbb{E}[f_i^-]=\frac{\mathbb{E}[|f_i|]}{2}\geq \frac{\varepsilon}{2} P[|f_i| \geq \varepsilon] \geq 3 \alpha \varepsilon.
\]
Denoting by $B_i=\{ f_i^- > \varepsilon \alpha\}$, H\"older's inequality implies
\[
\|f_i^-\|_{\infty} \mathbb{E} [1_{B_i}]\geq \mathbb{E}[f_i^-1_{B_i}]=\mathbb{E}[f_i^-]-\mathbb{E}[f_i^-1_{B_i^c}]\geq 3 \alpha \varepsilon- \varepsilon \alpha \geq 2 \alpha \varepsilon
\]
and the assertion follows since $\|f_i^-\|_{\infty}\leq \|f_i\|_{\infty}\leq 2 \varepsilon$.
\end{proof}

\section{The Kreps-Yan Theorem}

We follow here the proof as presented in~\cite{kab:97}:
\begin{theorem}\label{krepsyan}
Fix $p\in[1,\infty]$ and set $q$ conjugate to $p$. Suppose $C\subseteq L^p$ is a convex cone with $C\supseteq-L_{\geq 0}^p$ and $C\cap L_{\geq 0}^p=\{0\}$. If $C$ is closed in $\sigma(L^p,L^q)$, then there exists $ Q \sim P$ with $\frac{dQ}{dP}\in L^q(P)$ and $\mathbb{E}_Q[Y]\leq0 $ for all $ Y\in C$.
\end{theorem}
\begin{proof}
Any $x \in L_{\geq 0}^p\backslash\{0\}$ is disjoint from $C$, so we can by the Hahn-Banach-theorem strictly separate $x$ from $C$ by some $z_x\in L^q$. The cone property gives us $\mathbb{E}[z_xY]\leq0 $, for all $ Y\in C$ and $C \supset -L_{\geq 0}^p$ gives $z_x\geq0$. Strict separation implies $z_x\neq 0$, so that we can normalize to $\mathbb{E}[z_x]=1$.

We next form the family of sets $\{\Gamma_x:=\{z_x>0\}|x\in L_{\geq 0}^p\backslash\{0\}\}$. Then one can find a countable subfamily $(\Gamma_{x_i})_{i\in\mathbb{N}}$ with $P[\cup_i\Gamma_{x_i}]=1$. For suitably chosen weights $\gamma_i>0,\,i\in\mathbb{N}$, one gets that $Z:=\sum_{i=1}^{\infty}\gamma_iz_{x_i}$ is $Z>0$ almost surely with respect to $P$, $Z\in L^q$ and $\mathbb{E}[ZY]\leq0 $, for all $Y\in C$. Through normalization we get to $\mathbb{E}[Z]=1$, then $dQ:=ZdP$ does the job.
\end{proof}

\bibliographystyle{abbrv}

\bibliography{referencesFTAP140708}
\end{document}